\magnification=\magstephalf
\documentstyle{amsppt}
\NoBlackBoxes
\topmatter
\title
 On minimal poincar\'e  $4$--complexes
\endtitle
\author
Alberto Cavicchioli -- Friedrich Hegenbarth -- Du\v {s}an Repov\v {s} 
\endauthor
\subjclass 57 N 65, 57 R 67, 57 Q 10
\endsubjclass
\abstract We consider two types of minimal Poincar\'e $4$--complexes. One is defined with respect to the degree $1$--map order. This idea was already present in our previous  papers, and more systematically studied later by Hillman. The second type of minimal Poincar\'e  $4$--complexes were introduced by Hambleton, Kreck and Teichner. It is not based on an order relation. In the present paper we study  existence and uniqueness. 
\endabstract
\endtopmatter
\document

\subhead  1. Introduction \endsubhead

\medskip

Minimal objects are usually defined with respect to a partial order. We consider oriented Poincar\' e $4$--complexes (in short, $\operatorname{PD}_4$--complexes). If $X$ and $Y$ are two $\operatorname{PD}_4$--complexes, we define $X \succ Y$ if there is a degree $1$-map $f: X \to Y$ inducing an isomorphism on the fundamental groups. If also $Y \succ X$, well--known theorems imply that $f : X \to Y$ is a homotopy equivalence. So "$\succ $" defines a symmetric partial order on the set of homotopy types of $\operatorname{PD}_4$--complexes. A $\operatorname{PD}_4$--complex $P$ is said to be {\it minimal} for $X$ if $X \succ P$ and whenever $P \succ Q$, $Q$ is homotopy equivalent to $P$. We also consider special minimal objects called {\it strongly minimal}.  In this paper we study existence and uniqueness questions. It is an interesting question to calculate homotopy equivalences of $X$ relative a minimal $P$, that is, if $f: X\to P$ is as above, then calculate 
$$
\aligned
\operatorname{Aut}(X \succ P) =\{ h: X\to X : h \, \, & \text{homotopy equivalence such that} \, \, f \circ h \\
& \text{is homotopic to} \, \, f \}.
\endaligned
$$
Self--homotopy equivalences are studied by various authors (see [12] and references there). Pamuk's method can be used to calculate $\operatorname{Aut}(X \succ P)$.

\noindent
Constructions of minimal objects were indicated by Hegenbarth, Repov\u {s} and Spaggiari in [6] and more recently by Hillman in [8] and [9]. Degree $1$--maps can be constructed from $\Lambda$--submodules $G\subset H_2(X, \Lambda)$. More precisely, we have the following:

\proclaim{Proposition 1.1}  (cf. Proposition 2.4 below) Suppose $X$ is a Poincar\'e $4$--complex, and $G\subset
H_2(X,\Lambda)$  a  stably free $\Lambda$--submodule such that
the intersection form 
$\lambda_X$ restricted to $G$ is  nonsingular. Then one can  construct a Poincar\'e $4$--complex $Y$ and 
a degree $1$--map $f: X\to Y$. Moreover, there is an isomorphism
$$
\CD
K_2(f, \Lambda) = \operatorname{Ker} (H_2(X, \Lambda) @> f_{*} >> H_2(Y, \Lambda)) \cong G
\endCD
$$
and $\lambda_X$ restricted to $K_2(f, \Lambda)$ coincides with $\lambda_X$ on $G$ via this isomorphism.
\endproclaim

\proclaim{Corollary 1.2}  Given any   Poincar\'e $4$--complex $X$,
there exists a minimal Poincar\'e $4$--complex $P$ for $X$.
\endproclaim

The above proposition is useful to answer the following two basic questions about the minimal objects:

(1) Existence; and

(2) Uniqueness.

\medskip

A Poincar\'e $4$--complex $P$ is called {\it strongly minimal} for $\pi$ if the adjoint map $\overset \land \to{\lambda}_P : H_2 (P, \Lambda) \to  \operatorname{Hom}_{\Lambda} (H_2(P, \Lambda), \Lambda)$ of the intersection form $\lambda_P$ vanishes
[8]. Proposition 1.1 implies that $P$ is minimal.
The same questions arise if we consider the originally defined minimal objects in [5]. 

Existence of strongly minimal models $P$ is known only for few fundamental groups $\pi$ (see [5] and [8]). All these examples satisfy $H^3(B\pi, \Lambda) \cong 0$, hence $\operatorname{Hom}_{\Lambda} (H_2(P, \Lambda), \Lambda) \cong 0$ (see below). So all are "trivial" in the sense that $\lambda_P$ is zero because its adjoint  
$\overset \land \to{\lambda}_P : H_2(P, \Lambda) \to 
\operatorname{Hom}_{\Lambda} (H_2(P, \Lambda), \Lambda)
$ 
maps to the trivial $\Lambda$--module. 

An interesting question is therefore: {\sl Do there exist strongly minimal models} $P$ {\sl such that} $H^3(B\pi_1(P), \Lambda) \neq 0$?

We prove the following:

\proclaim{Theorem 1.3} Let $\pi$ be a finitely presented group such that 
$H^2(B\pi, \Lambda)$ is not a torsion group. Let $P$ and $P'$ be strongly minimal models for $\pi$.  Then $P$ and $P'$ are homotopy equivalent if the map $G : H_4(D, \Bbb Z)\to \operatorname{Hom}_{\Lambda} (H^2(D,\Lambda), \overline{H}_2(D, \Lambda))$ is injective, and if the $k$--invariants of $P$ and $P'$ correspond appropriately.
\endproclaim

Here $D$ is a $2$--stage Postnikov space and $G$ is defined via cap-products. Apart from the $k$--invariant, the injectivity of the map $G$ is an essential condition for uniqueness of strongly minimal models. In Section 4 we consider groups $\pi$ such that $B\pi$ is homotopy equivalent to a $2$--complex and prove that for any element of $\operatorname{Ker} G$ one can  construct a strongly minimal model. More precisely, we obtain:

\proclaim{Theorem 1.4} Suppose $B\pi$ is homotopy equivalent to a $2$--complex, and  
$\pi_2 = H^2(B\pi, \Lambda)$ is not a torsion group.  Then
$\operatorname{Ker} G\cong \Gamma(\pi_2)\otimes_{\Lambda} \Bbb Z$. Moreover, for any strongly minimal model $P$ and any $\xi \in \Gamma(\pi_2)$, it can be constructed another strongly minimal model $X$. 
\endproclaim

Examples are given by solvable Baumslag--Solitar groups (see [5]), or by surface fundamental groups. In Section 5 we construct non--homotopy equivalent strongly minimal models for these fundamental groups.

\bigskip

\subhead 2. Construction of degree $1$--maps \endsubhead

\bigskip

In this section we are going to prove Proposition 1.1 announced in Section 1. First we mention a result of Wall [14].

\proclaim{Lemma 2.1}  Let  $f : X\to Y$ be a
degree $1$--map between Poincar\'e $4$--complexes and suppose that $f_{*} : \pi_1(X)\to \pi_1(Y)$ is an isomorphism. 
Then $K_2(f, \Lambda) = \operatorname{Ker}(H_2(X, \Lambda)\to H_2(Y,\Lambda))$ is a stably $\Lambda$--free
submodule of $H_2(X, \Lambda)$ and $\lambda_X$ restricted to $K_2(f, \Lambda)$ is nonsingular. Also $K_2(f, \Lambda)\subset H_2(X, \Lambda)$ is a direct summand.
\endproclaim

\noindent
This section is devoted to proving a converse statement to Lemma 2.1.

\noindent
First we will show  Proposition 2.2. Before that, let us note that $\Lambda$ has an anti--involution that permits to switch from $\Lambda$--left to $\Lambda$--right modules, and to introduce compatible $\Lambda$--module structures on Hom--duals, etc.  We follow Wall's convention and consider $\Lambda$--right modules.

\proclaim{Proposition 2.2} Let $X$ be a Poincar\'e $4$--complex and $G\subset H_2(X, \Lambda)$ a $\Lambda$--free submodule so that $\lambda_X$ restricts to a nonsingular hermitian pairing on $G$. Then there exist a Poincar\'e $4$--complex $P$ and a degree $1$--map $f: X\to P$ such that $f_{*} : \pi_1(X)\to \pi_1(P)$ is an isomorphism and $K_2(f, \Lambda)\cong G$.  
\endproclaim

\demo{Proof} We recall that $\lambda_X$ is defined as the composite map
$$
\CD
H^{2}(X, \Lambda) \times H^2 (X, \Lambda) @> \cup >> H^4(X, \Lambda \otimes_{\Bbb Z} \Lambda)\cong H_{0}(X, \Lambda \otimes_{\Bbb Z} \Lambda) \cong \Bbb Z 
\otimes_{\Lambda} (\Lambda \otimes_{\Bbb Z} \Lambda) \\
@A \cong AA   @AA \cong A \\
H_{2}(X, \Lambda) \times H_2 (X, \Lambda) @> \lambda_X >>  \Lambda \cong \Lambda \otimes_{\Lambda} \Lambda 
\endCD
$$
and 
$$
\overset \land \to{\lambda}_X : H_2(X, \Lambda) \to \operatorname{Hom}_{\Lambda}(H_2(X, \Lambda), \Lambda)
$$
is the adjoint map of $\lambda_X$.

\noindent
To construct $P$, we consider a $\Lambda$--base $a_1$, \dots, $a_r$ of $G\subset  H_2(X, \Lambda)\cong \pi_2(X)$, and 
$$
\varphi_1, \dots, \varphi_r : \Bbb S^2 \to X
$$
representatives of $a_1$, \dots, $a_r$, respectively. Then $P$ is obtained from $X$ by adjoining $3$--cells along 
$\varphi_1$, \dots, $\varphi_r$. So $X\subset P$, and 
$$
H_p(P,X, \Lambda) \cong \cases G \quad p=3\\
0 \quad \text{otherwise} 
\endcases
\qquad \qquad  
H^p(P,X, \Lambda) \cong \cases G^{*} = \operatorname{Hom}_{\Lambda}(G, \Lambda) \qquad p=3\\
0 \qquad \qquad \qquad \qquad \text{otherwise}. 
\endcases
$$
Moreover, the sequence
$$
\CD
0 @>>> H_3(P, X, \Lambda) @> \partial_{*} >> H_2(X, \Lambda) @>>> H_2(P, \Lambda) @>>> 0
\endCD
$$
is exact because $\partial_{*} : H_3(P, G, \Lambda) \to G \subset H_2(X, \Lambda)$ is an isomorphism.

Note that there is a natural homomorphism
$$
\mu : H^2 (X, \Lambda) \to \operatorname{Hom}_{\Lambda}(H_2(X, \Lambda), \Lambda)
$$
such that the diagram
$$
\CD
H^2 (X, \Lambda) @> \mu >> \operatorname{Hom}_{\Lambda}(H_2(X, \Lambda), \Lambda)\\
@V  \cap [X] VV  @| \\
H_2 (X, \Lambda) @> \overset \land \to{\lambda}_{X} >> \operatorname{Hom}_{\Lambda}(H_2(X, \Lambda), \Lambda)
\endCD 
$$
commutes.
Let $[P] = f_{*}[X]$, where $f: X\subset P$ is the inclusion. Consider the diagrams
$$
\minCDarrowwidth{12pt}
\CD
0 @>>> H^2(P,\Lambda) @> f^{*} >> H^2 (X, \Lambda)  @> \delta^{*} >> H^3(P, X, \Lambda) @> \mu > \cong > 
\operatorname{Hom}_{\Lambda}  (H_3(P, X, \Lambda),\Lambda) = G^{*}\\
@. @V \cap [P] VV @VV \cap [X] V @. @AA \overset \land \to{\lambda}_{G} A \\
0 @<<< H_2(P, \Lambda) @< f_{*} << H_2(X, \Lambda) @< \partial_{*} << H_3(P, X, \Lambda) =G @= G @<<< 0    
\endCD
$$
and
$$
\CD
H^2 (X, \Lambda)  @> \delta^{*} >> H^3(P, X, \Lambda)  \\
@V \mu VV @V \cong V \mu V \\
\operatorname{Hom}_{\Lambda}  (H_2(X, \Lambda),\Lambda) @>>>
\operatorname{Hom}_{\Lambda}  (H_3(P, X, \Lambda),\Lambda) = G^{*}\\
@AA \overset \land \to{\lambda}_{X} A 
@AA \overset \land \to{\lambda}_{G} A \\
H_2(X, \Lambda) @< \partial_{*} << H_3(P, X, \Lambda) =G.     
\endCD
$$
Here $\overset \land \to{\lambda}_G = \overset \land \to{\lambda}_X |_{G}  $.
The left--hand square of the first diagram commutes. Combining  the right--hand square of the first diagram with the second diagram gives only
$$
\mu \circ \delta^{*} \circ (\cap [X])^{-1}\circ \partial_{*} = \overset \land \to{\lambda}_{G}.
$$
But this is sufficient to deduce  that $\cap [P] : H^2(P,\Lambda)\to H_2(P, \Lambda)$ is an isomorphism. It follows from the above short exact sequence that 
$$
\CD
f_{*} : H_3(X, \Lambda) @>> \cong > H_3(P, \Lambda) \qquad \qquad f^{*} : H^3(P, \Lambda) @>> \cong > H^3(X, \Lambda)
\endCD
$$
hence we obtain that 
$$
\CD
\cap  [P] : H^{*}(P, \Lambda) @>> \cong > H_{4 - *} (P, \Lambda)
\endCD
$$
for all $*$. The map $f$ is obviously of degree $1$.\qed
\enddemo

In the following we need another result of Wall about Poincar\'e complexes (see for instance [14]).

\proclaim{Lemma 2.3} Any Poincar\'e $4$-complex $X$ 
is homotopy equivalent to a CW--complex of the form $K\cup_{\varphi} D^4$, where $K$ is a $3$--complex and $\varphi : \Bbb S^3 \to K$ is an attaching map of the single $4$--cell $D^4$.
\endproclaim

Proposition 2.2 can be improved so that together with Lemma 2.1, we obtain the following:

\proclaim{Proposition 2.4} Let $X$ be a Poincar\'e $4$-complex. 
There exists a degree $1$--map $f : X \to Q$ if and only if there exists a stably free $\Lambda$--submodule $G\subset H_2(X, \Lambda)$ so that $\lambda_X$ restricts to a nonsingular hermitian form on $G$. In this case, $G \cong K_2(f, \Lambda)$.
\endproclaim

\demo{Proof} By Lemma 2.3 we can identify $X = K \cup_{\varphi} D^4$. The submodule $G$ is stably free, so $G\oplus H \cong \oplus_{1}^{\ell} \Lambda$, where $H$ is $\Lambda$--free. We may assume $H = \oplus_{1}^{2m} \Lambda$. Let $Z = X \# ( \#_{1}^{m} (\Bbb S^2 \times \Bbb S^2))$ be the Poincar\'e $4$-complex formed from $X$ by connected sum inside the $4$--cell with $\#_{1}^{m} (\Bbb S^2 \times \Bbb S^2)$.Then $G\oplus H \subset H_2(Z, \Lambda)$ and $\lambda_Z$ restricted to $H$ is the canonical hyperbolic form. If $a_1$, \dots, $a_{\ell}\in G\oplus H $  is a $\Lambda$--base, we attach $3$--cells to $Z$ along representatives $\varphi_1$, \dots, $\varphi_{\ell} : \Bbb S^2 \to X$ as in Proposition 2.2. We obtain a Poincar\'e $4$--complex $Q$ and a degree $1$--map $g : Z\to Q$ with $K_2(g, \Lambda) = G\oplus H$. We are going to show that $g$ factors over the collapsing map 
$$
c : Z = X \# ( \#_{1}^{m} (\Bbb S^2 \times \Bbb S^2)) \to X
$$
giving a degree $1$--map $f : X\to Q$. Note that
$$
 X \# ( \#_{1}^{m} (\Bbb S^2 \times \Bbb S^2)) \backslash  \text{4--cell} \simeq K \lor \{  \lor_{1}^{m} (\Bbb S^2 \lor \Bbb S^2)\}
$$
and the attaching map of the $4$--cell of $Z$ is of the following type
$$
a\oplus b \in \pi_3(K)\oplus [\pi_3(\lor_{1}^{m} (\Bbb S^2 \lor \Bbb S^2)) \otimes \Lambda] \subset \pi_3(Z \backslash 
(\text{4--cell}))
$$
where $a = [\varphi]$ and $b = [\psi] \otimes 1$ with $\psi : \Bbb S^3 \to \lor_{1}^{m} (\Bbb S^2 \lor \Bbb S^2)$ the attaching map of the $4$--cell of $\#_{1}^{m} (\Bbb S^2 \times \Bbb S^2)$. Obviously, $a\oplus b$ maps to zero in $\pi_3(Q)$. 

Now we apply Whitehead's $\Gamma$--functor to 
$$
\pi_2(Z) \cong \pi_2(K)\oplus H \cong \pi_2(Z \backslash (\text{4--cell})): \Gamma(\pi_2(Z)) \cong \Gamma(\pi_2(K)) \oplus \Gamma(H) \oplus \pi_2(K) \otimes H.
$$
The $\Gamma$--functor fits into certain Whitehead's exact sequence (see [1] and [15]) and by naturality one has the following diagram:
$$
\minCDarrowwidth{12pt}
\CD
0 @>>> \Gamma(\pi_2 (K)\oplus H) @>>> \pi_3 (Z\backslash (\text{4--cell})) @>>> H_3 (Z\backslash (\text{4--cell}), \Lambda) @>>> 0\\
@. @| @. @| \\
0 @>>>  \Gamma(\pi_2 (K)) \oplus \Gamma(H) \oplus \pi_2(K) \otimes H @>>> \pi_3 (Z\backslash (\text{4--cell})) @>>> 
H_3 (K, \Lambda) @>>> 0\\
@. @VVV @VVV @VVV \\
H_4(Q, \Lambda) @>>> \Gamma(\pi_2(Q)) @>>> \pi_3(Q) @>>> H_3(Q, \Lambda) @>>> 0.
\endCD
$$
Obviously, $b\in \Gamma(H) \subset \pi_3 (Z\backslash (\text{4--cell}))$, hence $b = \sum \lambda_{i j} [e_i, e_j]$, where 
the set
$\{e_1, \dots, e_{2m}\} \subset H$ is the standard base and $[ \cdot \, , \, \cdot ]$ denotes the Whitehead product. Now $H\subset G\oplus H\subset \pi_2(Z)$ maps to zero under $g_{*} :\pi_2(Z)\to \pi_2(Q)$, so $b \in \Gamma(\pi_2(K)\oplus H)$ maps to zero in $\Gamma(\pi_2(Q))$, hence it is zero in $\pi_3(Q)$. Because $a \oplus b$ is zero in $\pi_3(Q)$, also $a\in \pi_3(K)$ maps to zero under $\pi_3(K) \to \pi_3(Q)$. Therefore the inclusion map $K\subset Q$ extends to $f: X\to Q$, and $f$ induces a map 
$$
(X, K) \to (Q, Q\backslash (\text{4--cell})).
$$
We also have  
$$
g : (Z, Z\backslash (\text{4--cell}))\to
(Q, Q\backslash (\text{4--cell}))
$$
and a collapsing map
$$
c: (Z, Z\backslash (\text{4--cell}))\to (X, K).
$$
Since $Q$ is obtained from $Z$ by adding $3$--cells attached away from the $4$--cell, the following diagram commutes:
$$
\CD
H_4(Z, Z\backslash (\text{4--cell}), \Bbb Z)  @> c_{*} >> H_4(X, K, \Bbb Z) \\
@V g_{*} VV @VV f_{*} V \\
H_4(Q, Q\backslash (\text{4--cell}), \Bbb Z) @=
H_4(Q, Q\backslash (\text{4--cell}), \Bbb Z).
\endCD
$$
Because $c_{*}$ and $g_{*}$ map the fundamental class to the fundamental class, the degree of $f$ is $1$.\qed
\enddemo
\medskip

\noindent
{\bf Proof of Corollary 1.2}. We observe that for any degree $1$--map $f : X\to Y$ with $f_{*} : \pi_1(X) \underset \cong \to{\to} \pi_1(Y)$, one has 
$$
K_2(f, \Lambda) \otimes_{\Lambda} \Bbb Z = K_2(f, \Bbb Z) = \operatorname{Ker}(H_2(X, \Bbb Z)\to H_2(Y, \Bbb Z)),
$$
and that $H_2(X, \Bbb Z)$ is finitely generated. By Proposition 2.4 we can successively construct degree $1$--maps 
$$
\CD
X @> f >> Q, \qquad Q_1 @> f_1 >> Q_2, \quad \cdots
\endCD
$$
if we find nondegenerate stably free nontrivial submodules in $H_2(Q_k, \Lambda)$, and one has 
$$
K_2(f_k \circ \cdots \circ f_1 \circ f, \Lambda)\cong K_2(f_k, \Lambda) \oplus \cdots \oplus K_2(f_1, \Lambda) \oplus K_2(f, \Lambda) \subset H_2(X, \Lambda).
$$
Now 
$$
K_2(f_k \circ \cdots \circ f_1 \circ f, \Bbb Z)\cong K_2(f_k, \Bbb Z) \oplus \cdots \oplus K_2(f_1, \Bbb Z) 
\oplus K_2(f, \Bbb Z) \subset H_2(X, \Bbb Z)
$$
is finitely generated. Hence after certain $k$, we have 
$$
K_2(f_{k+1}, \Lambda) \otimes_{\Lambda} \Bbb Z = K_2(f_{k+1}, \Bbb Z) = \{0\}.
$$
Kaplansky's lemma (see remark below) implies $K_2(f_{k+1}, \Lambda) \cong 0$. Therefore $g = f_k \circ \cdots \circ f_1 \circ f : X\to Q_k$ is of degree $1$, and $Q_k$ is minimal. This completes the proof of Corollary 1.2.

\medskip

\noindent
{\bf Remark}. In [10], p.122, the following result is stated:

\medskip

\noindent
\proclaim{Lemma} Let $\Bbb F$ be a field of characteristic zero, and $\pi$ an arbitrary group. Let $A = \Bbb F[\pi]$ be the group algebra, and let $u$, $v\in M_n (A)$ be two $(n\times n)$ matrices such that the product $vu$ is the identity matrix $I_n$. Then $uv=I_n$.
\endproclaim

\noindent 
It has the following consequence (referred above as "Kaplansky's lemma"):

\medskip

\noindent
\proclaim{Corollary} If $K_2 (f, \Lambda) \otimes_{\Lambda} \Bbb Q \cong 0$, then $K_2 (f, \Lambda) \cong 0$.
\endproclaim

\demo{Proof} We know that $K_2 = K_2 (f, \Lambda)$ is stably free, i.e., $K_2 \oplus \Lambda^a \cong \Lambda^b$, where $a$ and $b$ are positive integers. Tensoring with $\Bbb Q$ implies that $a=b$. Let $h : K_2 \oplus \Lambda^a \to \Lambda^b$ be an isomorphism, and consider
$$
\CD
u = h \circ i : \Lambda^a @> \subset >>  K_2 \oplus \Lambda^a  @> h >> \Lambda^a 
\endCD
$$
and
$$
\CD
v = \operatorname{pr} \circ h^{-1} : \Lambda^a @> \subset >>  K_2 \oplus \Lambda^a  @> \operatorname{pr} >> \Lambda^a. 
\endCD
$$
Obviously $v\circ u = \operatorname{Id}$, hence $u\circ v = \operatorname{Id}$. This implies that $K_2 \subset \operatorname{Ker} (u\circ v) \cong 0$.
\qed\enddemo

\noindent
Note also that $K_2 \otimes_{\Lambda} \Bbb Q \cong 0$ is equivalent to $K_2 \otimes_{\Lambda} \Bbb Z\cong 0$.

\noindent
Of course, starting with $X$ one cannot in general assume that there is only one minimal $P$ and degree $1$--map $f: X\to P$ with $f_{*}: \pi_1(X) \underset \cong \to{\to} \pi_1(P)$.
\medskip

{\bf Problem 2.5}: Construct examples of $X$ which admit several minimal Poincar\'e $4$--complexes $P_i$ and degree $1$--maps $f_i : X\to 
P_i$ satisfying  $f_{i *}: \pi_1(X) \underset \cong \to{\to} \pi_1(P_i)$.

\bigskip

The next proposition completes the description of the correspondence between stably free $\Lambda$--modules with 
nondegenerate hermitian forms and degree $1$--maps of Poincar\'e $4$--complexes. But we have to assume that $\pi_1(X)$ does not contain elements of order $2$. 

\proclaim{Proposition 2.6} Let $X$ be a Poincar\'e $4$-complex and $G$ a stably free $\Lambda$--module with nondegenerate hermitian form. Then
there are a Poincar\'e $4$--complex $Y$ and 
a degree $1$--map $f : Y \to X$ such that  $K_2(f, \Lambda) \cong G$, $\lambda_Y$ restricted to
$K_2(f, \Lambda)$ coincides with $\lambda$ on $G$ under the isomorphism. Moreover, $f_{*} : \pi_1(Y)\to \pi_1(X)$ is an isomorphism.
\endproclaim

\demo{Proof} Let first $G$ be free of rank $m$. The proof procedes as in [7]. Here we begin with $Y' = X \# ( \#_{1}^{m} \Bbb CP^2)$ and the hermitian form $\lambda$, and continue as in Section 3 of [7] to construct $f: Y\to X$. If $G$ is stably free, that is, $G\oplus H \cong \Lambda^m$, where $H = \Bbb \Lambda^t$, we begin with $Y' = X \# ( \#_{1}^{m} \Bbb CP^2)$ and the hermitian form $\lambda' = \pmatrix \lambda & 0 \\ 0 & 1\endpmatrix $ on $G\oplus H$, and construct a degree $1$--map $f^{''} : Y^{''} \to X$ with $K_2(f^{''}, \Lambda) = G\oplus H$, and $\lambda_{Y^{''}}$ restricted to $K_2(f^{''},\Lambda)$ is equal to $\lambda'$. Now $H\subset H_2(Y^{''}, \Lambda)$ is $\Lambda$--free, and $\lambda_{Y^{''}}$ restricted to $H$ is non--singular. As in the proof of Proposition 2.4 we can construct a degree $1$--map $f: Y\to X$ with $K_2(f, \Lambda) \cong G$.
\qed
\enddemo

\subhead 3. A general result on the uniqueness of strongly minimal models
\endsubhead
\bigskip

Let $\pi$ be a finitely presented group. Suppose we are given strongly minimal Poincar\'e $4$--complexes $P$ and $P'$ with $\pi_1(P) \cong \pi \cong \pi_1(P')$. For simplicity, we denote $\pi'_1 =\pi_1(P')$, $\pi_1 = \pi_1(P)$, $\Lambda' = \Bbb Z [\pi'_1]$, and $\Lambda = \Bbb Z[\pi_1]$. Then we have 
$$
\CD
p^{*} : H^2(B\pi_1, \Lambda) @>> \cong > H^2(P, \Lambda)\\
p^{'*} : H^2(B\pi'_1, \Lambda') @>> \cong > H^2(P', \Lambda')
\endCD
$$
where $p: P\to B\pi_1$ and $p': P'\to B\pi'_1$ are the classifying maps. We denote by $\chi : D \to B\pi_1$ and 
$\chi' : D' \to B\pi'_1$ the $2$--stage Postnikov fibrations with fibres $K(\pi_2(P), 2)$ and $K(\pi_2(P'), 2)$, respectively. Spaces $D$ and $D'$ are obtained from $P$ and $P'$, respectively, by adding cells of dimension $\ge 4$ so that $\pi_q(D)\cong 0 \cong \pi_q(D')$ for every $q\ge 3$, and the inclusions $f: P\to D$ and $f': P'\to D'$ induce isomorphisms $f_{*} : \pi_i(P)\to \pi_i(D)$ and $f'_{*} : \pi_i(P')\to \pi_i(D')$, for every $i < 3$. We shall often write it as diagrams
$$
\CD
P @> f >> D \\
@| @VV \chi V \\
P @>p >> B\pi_1
\endCD \qquad \qquad  
\CD
P' @> f' >> D' \\
@| @VV \chi' V \\
P' @>p' >> B\pi'_1.
\endCD 
$$

We choose an isomorphism $\alpha : \pi_1 \to \pi'_1$. It determines an isomorphism $\Lambda \to \Lambda'$ of rings. For the sake of simplicity we shall identify $\Lambda'$ with $\Lambda$ via this isomorphism when we use it as coefficients in (co)homology groups. 
We define 
$$
\beta : H_2(P,\Lambda) \to H_2(P', \Lambda)
$$
by the following diagram
$$
\CD
H^2(B\pi_1, \Lambda) @> p^{*} > \cong > H^2(P, \Lambda) @> \cap [P] > \cong > H_2(P, \Lambda)\\
@A (B\alpha)^{*}  AA     @.   @VV \beta V \\
H^2(B\pi'_1, \Lambda) @> p^{'*} > \cong > H^2(P', \Lambda) @> \cap [P'] > \cong > H_2(P', \Lambda).
\endCD
\tag 1
$$
The next diagram explains the compatibility of the $k$--invariants $k^3_P \in H^3(B\pi_1, \pi_2(P))$ and 
$k^3_{P'} \in H^3(B\pi'_1, \pi_2(P'))$
$$
\CD
\operatorname{Hom}_{\Lambda}(H_2(P, \Lambda), H_2(P, \Lambda)) @>>> H^3(B\pi_1, H_2(P, \Lambda)) \\
@VV \beta_{\#} V @VV \beta_{\#} V \\
\operatorname{Hom}_{\Lambda}(H_2(P, \Lambda), H_2(P', \Lambda)) @>>> H^3(B\pi_1, H_2(P', \Lambda)) \\
@AA \beta^{\#} A @AA (B\alpha)^{*} A \\
\operatorname{Hom}_{\Lambda}(H_2(P', \Lambda), H_2(P', \Lambda)) @>>> H^3(B\pi'_1, H_2(P', \Lambda)) 
\endCD
\tag 2
$$
where the top (resp. bottom) horizontal map sends $\operatorname{Id}$ into $k^3_P$ (resp. $k^3_{P'}$), and on the left (resp. right) vertical side we have $\beta_{\#} (\operatorname{Id}) = \beta = \beta^{\#} (\operatorname{Id})$ (resp. 
$\beta_{\#} (k^3_P) =  (B\alpha)^{*} (k^3_{P'})$). Therefore there is a homotopy equivalence $h : D\to D'$ such that
the diagram
$$
\CD
D @> h >> D'\\
@VV \chi V @VVV \\
B\pi_1 @> B\alpha >> B\pi'_1
\endCD
$$
commutes (up to homotopy). Furthermore, Diagram (1) can be completed to the following diagram 
$$
\minCDarrowwidth{12pt}
\CD
H^2(D, \Lambda) @< \chi^{*} < \cong <  H^2(B\pi_1, \Lambda) @> p^{*} >> H^2(P, \Lambda) @> \cap [P] >> H_2(P, \Lambda) @> f_{*} >> H_2(D, \Lambda) \\
@A h^{*} AA @AA (B\alpha)^{*} A @. @V  \beta  VV  @VV h_{*} V \\
H^2(D', \Lambda) @< \chi^{'*} < \cong <  H^2(B\pi'_1, \Lambda) @> p^{'*} >> H^2(P', \Lambda) @> \cap [P'] >> H_2(P', \Lambda) @> f'_{*} >> H_2(D', \Lambda)
\endCD
\tag 3
$$
where
$$
\CD
H^2(D, \Lambda) @> f^{*} >> H^2(P, \Lambda) \\
@| @AA p^{*} A \\
H^2(D, \Lambda) @< \chi^{*} < \cong < H^2(B\pi_1, \Lambda)
\endCD
\qquad \qquad 
\CD
H^2(D', \Lambda) @> f^{'*} >> H^2(P', \Lambda) \\
@| @AA p^{'*} A \\
H^2(D', \Lambda) @< \chi^{'*} < \cong < H^2(B\pi'_1, \Lambda).
\endCD
$$
Note that all the maps are $\Lambda$--isomorphisms. 

At this point it is convenient to introduce the map
$$
G : H_4(D, \Bbb Z) \to \operatorname{Hom}_{\Lambda} (H^2(D, \Lambda), \overline{H}_2(D, \Lambda))
$$
using the equivariant cap-product construction, and similarly $G'$ for $D'$. From Diagram (3) we summarize as follows:

\proclaim{Corollary 3.1} Diagram (3) commutes, and the composed horizontal homomorphisms (from left to right) are $G(f_{*}[P])$ and $G'(f'_{*}[P'])$.
\endproclaim

We again invoke  Wall's theorem (Lemma 2.3) and identify 
$$
P = K \cup_{\varphi} D^4 \qquad \qquad  P' = K' \cup_{\varphi'}  D^{'4}
$$
where $K$ and $K'$ are $3$--complexes, and $\varphi : \Bbb S^3 \to K$ and $\varphi' : \Bbb S^3 \to K'$ are the attaching maps of the $4$--cells $D^4$ and $D^{'4}$, respectively. Hence $(D,K)$ and $D', K')$ are relative CW--complexes with cells in dimensions $k\ge 4$, that is, $D^{(3)} = K$ and ${D'}^{(3)} = K'$. Approximate $h : D\to D'$ by a cellular map (again denoted by $h$). Then 
$$
h^{(3)} = h|_{K} : K\to K'
$$
and
$$
\CD
D @> h >> D'\\
@A i AA @AA i'  A\\
K @> h^{(3)} >> K' 
\endCD
$$   
commutes, where $i: K\subset D$ and $i' : K' \subset D'$ are the inclusion maps.
 
\proclaim{Proposition  3.2} 

(a) $h^{(3)} : K\to K'$ extends to $\phi : P\to P'$ if $h_{*} f_{*} [P] = \ell f'_{*} [P'] \in H_4(D', \Bbb Z)$ for some $\ell \in \Bbb Z$; and

(b) If $f'_{*} : H_4(P', \Bbb Z) \to H_4(D', \Bbb Z)$ is  injective and $\ell = \pm 1$, then $\phi$ is of degree $\pm 1$, hence it is a homotopy equivalence.
\endproclaim

\demo{Proof} 
(a) The obstruction to extending $h^{(3)}$ belongs to 
$$
\aligned
H^4(P, \pi_3(P')) & \cong H_{0}(P, \pi_3(P')) \cong \Bbb Z \otimes_{\Lambda} \pi_3(P') \cong \Bbb Z \otimes_{\Lambda} \pi_4(D', P') \\
& \cong \Bbb Z \otimes_{\Lambda} H_4(D', P', \Lambda) = H_4(D', P', \Bbb Z))
\endaligned
$$
(one applies among others: $\pi_3(D')= \pi_3(D) = 0$ and the Hurewicz theorem). The obstruction in $\Bbb Z \otimes_{\Lambda} \pi_3(P')$ is given by the image of $[h^{(3)} \circ \varphi] \in \pi_3(K')$ under the composite map
$$
\CD
\pi_3(K') @>>> \pi_3(P') @>>> \pi_3(P') \otimes_{\Lambda} \Bbb Z.
\endCD
$$
The obstruction in $H_4(D', P', \Bbb Z)$ is given by the induced map of the composition
$$
\CD
(D^4, \Bbb S^3) @> \varphi >> (P,K) \subset (D,K) @> h >> (D', K') \subset (D', P')
\endCD
$$
hence it is the image of $[P] \in H_4(P, \Bbb Z)$ under the composition on the bottom horizontal row in the following diagram
$$
\CD
H_4(P,K, \Bbb Z) @>>> H_4(D,K,\Bbb Z) @> h_{*} >> H_4(D', K', \Bbb Z) @= H_4(D', K', \Bbb Z)\\
@A \cong AA @AAA @AAA @VVV \\
H_4(P, \Bbb Z) @> f_{*} >> H_4(D, \Bbb Z) @> h_{*} >> H_4(D', \Bbb Z) @>>> H_4(D', P', \Bbb Z).
\endCD
\tag 4
$$
Hence the obstruction vanishes if and only if $h_{*}f_{*} [P] = \ell f'_{*} [P']$ for some $\ell \in \Bbb Z$.
\medskip

(b) If $\phi : P\to P'$ exists, then it is such that the diagram
$$
\CD
H_4(P, \Bbb Z) @> \phi_{*} >> H_4(P', \Bbb Z) @> f'_{*} >> H_4(D', \Bbb Z)\\
@VV \cong V  @VV \cong V @| \\
H_4(P, K, \Bbb Z) @> \phi_{*} >> H_4(P', K', \Bbb Z) @. H_4(D', \Bbb Z)\\
@VV f_{*} V @VV f'_{*} V @| \\
H_4(D, K, \Bbb Z) @> h_{*} >> H_4(D', K', \Bbb Z) @<<< H_4(D', \Bbb Z)
\endCD
$$ 
commutes. Hence $f'_{*} \phi_{*} [P] = h_{*} f_{*} [P] = \pm f'_{*} [P']$ implies $\phi_{*} [P] = \pm [P']$ since $f'_{*}$ is injective. Using the Poincar\'e duality one obtains
$$
\CD
\phi_{*} : H_{*}(P, \Lambda) @>> \cong > H_{*}(P', \Lambda).
\endCD
$$
Because $\phi_{*} : \pi_1(P) \to \pi_1(P')$ is an isomorphism, the map $\phi : P\to P'$ is a homotopy equivalence by the Hurewicz--Whitehead theorem.
\qed
\enddemo

\bigskip

\noindent
{\bf Proof of Theorem 1.3}. We have a commutative diagram (up to homotopy)
$$
\CD
D @= D @> h >> D'  @= D' \\
@A f AA @VVV @VVV @AA f' A \\
P @> p >> B\pi_1 @> B\alpha >> B\pi'_1 @< p' << P'
\endCD
$$
where $h: D\to D'$ is a homotopy equivalence. Consider the diagram
$$
\CD
H_4(D,\Bbb Z) @> G >>
\operatorname{Hom}_{\Lambda} (H^2(D,\Lambda), \overline{H}_2(D, \Lambda))\\
@V h_{*} VV @VV T V \\
H_4(D',\Bbb Z) @> G' >>
\operatorname{Hom}_{\Lambda} (H^2(D',\Lambda), \overline{H}_2(D', \Lambda))
\endCD
\tag 5
$$
where $\cap z$ is the cap product with $z\in H_4(D, \Bbb Z)$. Similarly, $\cap'$. The map $T$ is defined by $T(\xi) = h_{*} \circ \xi \circ h^{*}$. Note that $T$ is an isomorphism.\qed

\proclaim{Lemma 3.3} Diagram (5) commutes.
\endproclaim

\demo{Proof}  Given $x\in H_4(D, \Bbb Z)$ and $u' \in H^2(D', \Bbb Z)$, then we have
$$
TG(x)(u') = h_{*} (h^{*}(u') \cap x) = u' \cap h_{*}(x) = G'h_{*}(x)
$$
as required. 
\qed\enddemo

\medskip

\noindent
Now consider the diagram
$$
\CD
H_4(P, \Bbb Z) @> f_{*} >> H_4(D, \Bbb Z) @> G >> \operatorname{Hom}_{\Lambda} (H^2(D, \Lambda), \overline{H}_2(D, \Lambda))\\
@. @VV h_{*} V @VV T V \\
H_4(P', \Bbb Z) @> f'_{*} >> H_4(D', \Bbb Z) @> G' >> \operatorname{Hom}_{\Lambda} (H^2(D', \Lambda), \overline{H}_2(D', \Lambda)).
\endCD
$$
It follows from Corollary 3.1 that
$$
TG f_{*} [P] = G' f'_{*} [P'],
$$
and from $TG = G' h_{*}$ we get $G' h_{*} f_{*} [P] = G' f'_{*} [P']$, hence $h_{*} f_{*} [P] = f'_{*} [P']$. So Proposition 3.2 (a) holds with $\ell =1$.

A similar diagram as (5) holds for the space $P'$:

$$
\CD
H_4(P',\Bbb Z) @> G'' >>
\operatorname{Hom}_{\Lambda} (H^2(P',\Lambda), \overline{H}_2(P', \Lambda)) \cong 
\operatorname{Hom}_{\Lambda} (\overline{H}_2(P',\Lambda), \overline{H}_2(P', \Lambda))\\
@V f'_{*} VV @VV T V \\
H_4(D',\Bbb Z) @> G' >>
\operatorname{Hom}_{\Lambda} (H^2(D',\Lambda), \overline{H}_2(D', \Lambda))
\endCD
$$
with $T(\xi) = f_{*} \circ \xi \circ f^{*}$. Since $T$ is an isomorphism, $f'_{*}$ is injective if and only if the map 
$G''$ is injective. Now observe that under the maps
the generator $[P']$ goes to $\operatorname{Id}$. The upper right isomorphism is induced by Poincar\'e duality.
Hence $G''$ is injective if and only if $\operatorname{Id}$ is not of finite order. Now $H_2(P', \Lambda) \cong H^2(B\pi'_1, \Lambda) \cong H^2(B\pi_1, \Lambda)$. 
The claim now follows from  Proposition 3.2(b).
 
\medskip

\subhead 4. Construction of strongly minimal models
\endsubhead
\medskip

The principal examples of fundamental groups $\pi$ admitting a strongly minimal model $P$ are discussed in [5]. These are groups of geometric dimension equal to $2$, i.e., $B\pi$ is a $2$--dimensional aspherical complex. It is easy to see that the boundary of a regular neighbourhood $N$ of an embedding $B\pi \subset \Bbb R^5$ is a strongly minimal model for $\pi$ (see [5]). Here we show that the map $G$ is not injective, hence we cannot expect uniqueness up to homotopy equivalence. In fact, we are going to classify all strongly minimal models fixing $\pi$ by elements of the kernel of $G$. Note that all $k$--invariants vanish since $B\pi$ is a $2$--complex. We assume $H_4(P, \Lambda) \cong 0$, i.e., that $\pi$ is infinite (which holds for the known examples).
\medskip

(4.1) {\sl Computation of} $\operatorname{Ker} G$.
\medskip

We fix $\pi$ as above, and for conveniency also one strongly minimal model $P$, say $P =\partial N$. We have the $2$--stage Postnikov system
$$
\CD
D @> \chi >> B\pi \\
@A f AA @AA p A \\
P @= P
\endCD
$$

\proclaim{Lemma 4.1} There is an exact sequence
$$
\CD
0 @>>> \Gamma(\pi_2)\otimes_{\Lambda} \Bbb Z @>>> H_4(D,\Bbb Z) @>>> H_2(B\pi, H_2(D, \Lambda)) @>>> 0
\endCD
$$
where $\pi_2 =\pi_2(P)\cong \pi_2(D)$
\endproclaim

\demo{Proof}  This follows from the spectral sequence
$$
E^2_{p q} = H_p(B\pi, H_q(D, \Lambda)) \underset p+q = n \to{\Longrightarrow} H_n(D, \Bbb Z).
$$
Taking $n=4$ we have $E^2_{p q} = E^{\infty}_{p q} = [F_p H_4(D, \Bbb Z)]/[F_{p-1} H_4(D, \Bbb Z)]$ with filtration
$$
0\cong F_{-1} H_4 \subset F_{0} H_4 \subset F_{1} H_4 \subset F_{2} H_4 \subset F_{3} H_4 \subset F_{4} H_4(D, \Bbb Z) = H_4(D, \Bbb Z).
$$
The result follows since $E^2_{2 2} = H_2(B\pi, H_2(D, \Lambda))$, $E^2_{0 4} = H_{0} (B\pi, H_4(D, \Lambda)) = H_4(D, \Lambda)\otimes_{\Lambda} \Bbb Z$, and $E^2_{p q} \cong 0$ else for $p+q =4$.
\qed\enddemo

\medskip

\noindent
{\bf Remark.} Similarly one gets the exact sequence 
$$
\CD
0 @>>> H_1(P, H_3(P, \Lambda)) @>>> H_4(P, \Bbb Z) @>>> H_2(B\pi, H_2(P, \Lambda)) @>>> 0 .
\endCD
$$
In particular, $H_2(B\pi, H_2(D, \Lambda))$ is a quotient of $\Bbb Z$ because $H_2(D, \Lambda)\cong H_2(P, \Lambda)$ and $H_4(P, \Bbb Z) \cong \Bbb Z$.
\medskip

\proclaim{Lemma 4.2} The kernel of 
$$
G : H_4(D, \Bbb Z) \to \operatorname{Hom}_{\Lambda - \Lambda}(H^2(D, \Lambda), \overline{H}_2(D,\Lambda))
$$
is $\Gamma(\pi_2)\otimes_{\Lambda} \Bbb Z$. 
\endproclaim

\demo{Proof}  The map $\chi^{*} : H^2(B\pi, \Lambda) \to H^2(D, \Lambda)$ is an isomorphism, and $H^2(B\pi, \Lambda) \cong [\operatorname{Hom}_{\Lambda}(C_2(\widetilde {B\pi}), \Lambda)]/[\operatorname{Im} \delta^1]$, where
$$
\delta^1 : \operatorname{Hom}_{\Lambda} (C_1 ( \widetilde {B\pi}), \Lambda) \to \operatorname{Hom}_{\Lambda} (C_2 (\widetilde {B\pi}), \Lambda)
$$
is the coboundary map. The composition 
$$
\CD
\operatorname{Hom}_{\Lambda - \Lambda} (H^2(B\pi, \Lambda), \overline{H}_2(D, \Lambda)) @>>> \operatorname{Hom}_{\Lambda} (\operatorname{Hom}_{\Lambda} (C_2 (\widetilde {B\pi}), \Lambda), \overline{H}_2(D, \Lambda))\\
@A \cong AA \\
\operatorname{Hom}_{\Lambda}(H^2(D, \Lambda), H_2(D, \Lambda)) 
\endCD
$$
is obviously injective. Because $C_2(\widetilde{B\pi})$ is $\Lambda$--free there is a canonical isomorphism
$$
\operatorname{Hom}_{\Lambda - \Lambda} (\operatorname{Hom}_{\Lambda} (C_2 (\widetilde {B\pi}), \Lambda), \overline{H}_2(D, \Lambda)) \cong C_2 (\widetilde {B\pi}) \otimes_{\Lambda} H_2(D, \Lambda).
$$
Composing all these maps gives an injective map
$$
\operatorname{Hom}_{\Lambda - \Lambda} (H^2(D, \Lambda), \overline{H}_2(D, \Lambda)) \to C_2 (\widetilde {B\pi}) \otimes_{\Lambda} H_2(D, \Lambda).
$$
The composition with $G$ gives a map $H_4(D, \Bbb Z) \to C_2( \widetilde {B\pi}) \otimes_{\Lambda} H_2(D, \Lambda)$ with image the $2$--cycle subgroup of the complex $C_{*}(\widetilde {B\pi}) \otimes_{\Lambda} H_2(D, \Lambda)$, i.e., $H_2(B\pi, H_2(D, \Lambda))$. This is the map $H_4(D, \Bbb Z) \to H_2(B\pi, H_2(D, \Lambda))$ of Lemma 4.1. In other words, we have the following commutative diagram
$$
\CD
\operatorname{Hom}_{\Lambda - \Lambda} (H^2(D, \Lambda), \overline{H}_2(D, \Lambda)) @>>> H_2(B\pi, H_2(D, \Lambda)) \\
@A G AA @AAA \\
H_4(D, \Bbb Z) @= H_4(D, \Bbb Z) 
\endCD
$$ 
where the horizontal map is injective. The result now follows from Lemma 4.1.
\qed\enddemo
\medskip

{\sl Supplement to Lemma 4.1}. If $P$ and $P'$ are two strongly minimal models for $\pi$, let
$$
\CD
P @> f >> D \\
@| @VV \chi V \\
P @>> p > B\pi
\endCD
\qquad \qquad
\CD
P @> f' >> D' \\
@| @VV \chi'  V\\
P' @>> p' > B\pi
\endCD
$$
be the two associated $2$--stage Postnikov systems. Let $h : D\to D'$ be the homotopy equivalence constructed in Section 3. Then 
the diagram 
$$
\CD
H_4(D, \Bbb Z) @>>> H_2(B\pi, H_2(D, \Lambda)) \\
@V h_{*} VV @VVV \\
H_4(D', \Bbb Z) @>>> H_2(B\pi, H_2(D', \Lambda))
\endCD
$$
commutes. The right vertical map is induced by $h_{*} : H_2(D, \Lambda) \to H_2(D', \Lambda)$.
\medskip

(4.2) {\sl Construction of strongly minimal models}

\medskip

We choose a strongly minimal  model $P$ for $\pi$. By Wall's theorem [13], $P$ is homotopy equivalent to $K\cup_{\varphi_1} D^4$, where $K$ is a $3$--complex, and $\varphi_1 : \Bbb S^3 \to K$ is the attaching map of the only $4$--cell. This representation is unique, i.e., given a homotopy equivalence 
$$
\CD
K_1 \cup_{\varphi_1} D^4 @> h >> K_2 \cup_{\varphi_2} D^4
\endCD
$$
then there is a homotopy equivalence of pairs $(K_1, \varphi_1(\Bbb S^3)) \to (K_2, \varphi_2 (\Bbb S^3))$ (see [13], p.222). We simply write $P = K\cup_{\varphi_1} D^4$, and change the attaching map $[\varphi_1] \in \pi_3(K)$ by an element $[\varphi] \in \Gamma(\pi_2)$, i.e., $[\varphi] \in \Gamma(\pi_2) =\operatorname{Im} (\pi_3(K^{(2)}) \to \pi_3(K))$, and we consider $X = K \cup_{\varphi_2} D^4$, where $\varphi_2 = \varphi_1 + \varphi $ and $\varphi : \Bbb S^3 \to K^{(2)}$. Let $q : X\to B\pi$ be the classifying map. It follows that $q^{*} : H^2(B\pi, \Lambda)\to H^2(X, \Lambda)$ is an isomorphism. If $X$ is a Poincar\'e $4$--complex, then $X$ is a strongly minimal model for $\pi$.

\medskip

(4.3) {\sl Proof of the Poincar\' e duality}
\medskip

(I) We have an isomorphism $\pi_4(X, K) \to H_4(X, K, \Lambda)\cong \Lambda$.  Let us consider the diagram of Whitehead's sequences
$$
\CD
@. @. 0 @. 0 \\
@. @. @VVV @VVV \\
@. @. \pi_4(X, K) @> \cong >> H_4(X, K, \Lambda) \\
@. @. @VVV @VVV \\
0  @>>> \Gamma(\pi_2) @>>> \pi_3(K) @>>> H_3(K, \Lambda) @>>> 0 \\
@. @| @VVV @VVV  \\
@. \Gamma(\pi_2) @>>> \pi_3(X) @>>> H_3(X, \Lambda) @>>> 0 \\
@. @. @VVV @VVV \\
@. @. 0 @. 0 \\
\endCD
$$  
One has a similar diagram if we replace $X$ by $P$. Under the Hurewicz map, $[\varphi_1]$ and $[\varphi_2]$ go to the same element in $H_3(K, \Lambda)$ which coincides with the images of the generators of $H_4(P, K, \Lambda)$ resp. $H_4(X, K, \Lambda)$ under the connecting homomorphism, hence $H_3(X, \Lambda) \cong H_3(P, \Lambda)$. Moreover, this gives us the following:

\medskip

\proclaim{Lemma 4.3} $H_4(X, \Bbb Z)\cong \Bbb Z$
\endproclaim

\demo{Proof} 
Tensoring with $\otimes_{\Lambda} \Bbb Z$ the upper part of the above diagram gives
$$
\CD
@. \pi_4(X, K)\otimes_{\Lambda} \Bbb Z @>> \cong > H_4(X, K, \Lambda)\otimes_{\Lambda} \Bbb Z @>> \cong > H_4(X, K, \Bbb Z)\\
@. @VVV @VVV @VVV \\
\Gamma(\pi_2) \otimes_{\Lambda} \Bbb Z @>>> \pi_3(K)\otimes_{\Lambda} \Bbb Z @>>> H_3(K, \Lambda) \otimes_{\Lambda} \Bbb Z @>>> H_3(K, \Bbb Z)
\endCD
$$ 
Similarly, for $X$ replaced by $P$. (We don't claim the exactness of the lower row). Now $H_4(P, K, \Bbb Z)\to H_3(K, \Bbb Z)$ is the zero map. By the argument above, $[\varphi_1]\otimes_{\Lambda} 1$ and $[\varphi_2]\otimes_{\Lambda} 1$ map to the same element in $H_3(K, \Lambda) \otimes_{\Lambda} \Bbb Z$, hence the generators of $H_4(X, K, \Bbb Z)$ resp. $H_4(P, K, \Bbb Z)$ map to the same element in $H_3(K, \Bbb Z)$ under the connecting homomorphisms. Thus $H_4(X, K, \Bbb Z) \to H_3(K, \Bbb Z)$ is the zero map. Therefore there is an isomorphism $H_4(X, \Bbb Z) \to H_4(X, K, \Bbb Z) \cong \Bbb Z$.
\qed\enddemo

Let $[X]\in H_4(X, \Bbb Z)$ be a generator. We have to study 
$$
\cap [X] : H^p(X, \Lambda)\to H_{4-p}(X, \Lambda).
$$
To examine the cases $p=1$ and $p=3$, we introduce an auxiliary space $Y = K \cup_{\varphi_1, \varphi} \{D^4, D^4 \}$, obtained from $K$ by attaching two $4$--cells with attaching maps $\varphi_1$ and $\varphi$. Note that $Y = P\cup_{\varphi} D^4$.
\medskip

(II) {\sl Case} $p=1$
\medskip

Let $i : P\to Y$ be the inclusion, and $j : X\to Y$ be the map induced by $K\subset Y$ and 
$$
\CD
\varphi_2 = \varphi_1 + \varphi : \Bbb S^3 @>>> \Bbb S^3 \lor \Bbb S^3 @> \varphi_1 \lor \varphi >> K .
\endCD
$$ 
We have the following maps of pairs
$$
\CD
(D^4, \Bbb S^3) @> \bar i \circ \bar \varphi_1 >> (Y, K) \\
@V \bar \varphi_1 VV @AA \bar i A \\
(P, K) @= (P, K)
\endCD
\qquad \qquad
\CD
(D^4, \Bbb S^3) @> \bar j \circ \bar \varphi_2 >> (Y, K) \\
@V \bar \varphi_2 VV @AA \bar j A \\
(X, K) @= (X, K)
\endCD
$$
and $\bar \varphi : (D^4, \Bbb S^3) \to (Y, K)$. Obviously, $\bar \varphi_2 = \bar \varphi_1 + \bar \varphi : (D^4, \Bbb S^3) \to (Y, K)$ is the $4$--cell $[\varphi_1]$ "slided" over $[\varphi]$. Since $[\bar \varphi] \in \Gamma(\pi_2)$, $\bar \varphi$ factors as follows
$$
\CD
(D^4, \Bbb S^3) @> \bar k \circ \bar \varphi >> (Y, K) \\
@V \bar \varphi VV @AA  \bar k A \\
(K^{(2)}\cup_{\varphi} D^4, K^{(2)} ) @= (K^{(2)}\cup_{\varphi} D^4, K^{(2)} )
\endCD
$$

From this one sees that $\bar j_{*} [\bar \varphi_2] - \bar i_{*} [\bar \varphi_1]$ belongs to 
$$
\operatorname{Im} (H_4(K^{(2)} \cup_{\varphi} D^4, K^{(2)}) \to H_4(Y, K)).
$$
The diagram
$$
\CD
H_4(X) @> j_{*} >> H_4(Y) @< i_{*} << H_4(P) \\
@V \cong VV @VVV @VV \cong V \\
H_4(X, K) @> \bar j_{*} >> H_4(Y, K) @<  \bar i_{*} << H_4(P, K)
\endCD
$$
as well as injectivity of $H_4(Y)\to H_4(Y, K)$ and the isomorphism 
$$
H_4(K^{(2)} \cup_{\varphi} D^4) \to 
H_4(K^{(2)} \cup_{\varphi} D^4, K^{(2)}) 
$$ 
prove the following

\proclaim{Lemma 4.4} $j_{*} [X] - i_{*} [P]$ belongs to $\operatorname{Im} (H_4(K^{(2)} \cup_{\varphi} D^4) \to H_4(Y))$. 
\endproclaim

\medskip

\proclaim{Corollary 4.5} Taking cap products with $i_{*} [P]$ and $j_{*} [X]$ : $H^1(Y, \Lambda) \to H_3(Y, \Lambda)$ gives the same map.
\endproclaim

\demo{Proof} Let $\theta \in H_4(K^{(2)} \cup_{\varphi} D^4)$ map to $j_{*} [X] - i_{*} [P]$. Then the diagram
$$
\CD
H^1(Y, \Lambda) @> \cap j_{*} [X] - \cap i_{*} [P] >> H_3(Y, \Lambda) \\
@V \cong VV @AAA \\
H^1(K^{(2)} \cup_{\varphi} D^4, \Lambda) @> \cap \theta >> H_3(K^{(2)} \cup_{\varphi} D^4, \Lambda) \cong 0
\endCD
$$
commutes.
\qed\enddemo

\proclaim{Lemma 4.6}  $i_{*} : H_3(P, \Lambda)\to H_3(Y, \Lambda)$ is an isomorphism.
\endproclaim

\demo{Proof} Since $Y = P\cup_{\varphi} D^4$, $i_{*}$ is surjective. Let us consider the diagram
$$
\CD
H_4(K^{(2)} \cup_{\varphi} D^4, K^{(2)}, \Lambda) @>>> H_4(Y, P, \Lambda) @>>> H_3(P, \Lambda)\\
@A \cong AA @A \cong AA @AAA \\
H_4(K^{(2)} \cup_{\varphi} D^4,\Lambda) @>> \cong > 
H_4(K^{(2)} \cup_{\varphi} D^4, K^{(2)}, \Lambda) @>>> 
H_3(K^{(2)} \cup_{\varphi} D^4, \Lambda) 
\endCD
$$
which shows that $H_4(Y, P, \Lambda) \to H_3(P, \Lambda)$ is the zero map.
\qed\enddemo

\proclaim{Lemma 4.7}  $j_{*} : H_3(X, \Lambda)\to H_3(Y, \Lambda)$ is an isomorphism.
\endproclaim

\demo{Proof} The map $j_{*}$ is surjective because $Y^{(3)} =K = X^{(3)}$. We identify $H_4(Y, K, \Lambda) \equiv \Lambda \oplus \Lambda$ according to the diagram
$$
\minCDarrowwidth{12pt}
\CD
@. 0 @. 0 @. 0 \\
@. @VVV @VVV @AAA \\
H_4(D^4, \Bbb S^3, \Lambda) @> \bar \varphi_{1 *} > \cong > H_4(P, K, \Lambda) @. H_4(Y, \Lambda) @. H_4(Y, X, \Lambda)\\
@. @V \bar i_{*} VV @VVV @AAA  \\
@. H_4(Y, K, \Lambda) @= H_4(Y, K, \Lambda) @= H_4(Y, K, \Lambda) @< \bar k_{*} \bar \varphi_{*} << H_4(D^4, \Bbb S^3, \Lambda) \\
@. @A \bar j_{*} AA @. @VVV @VV \cong V \\
H_4(D^4, \Bbb S^3, \Lambda) @> \bar \varphi_{1 *} + \bar \varphi_{*} > \cong > H_4(X, K, \Lambda) @. @. H_4(Y, P, \Lambda) @= H_4(Y, P, \Lambda)  \\
@. @. @. @VVV \\
@. @. @. 0  
\endCD
$$
where $\bar i_{*} [\bar \varphi_1 ]= (1,0) \in \Lambda \oplus \Lambda$ and $\bar k_{*} [\bar \varphi]= (0, 1) \in \Lambda \oplus \Lambda$. The map $\bar k_{*} \bar \varphi_{*}$ defines a splitting of $H_4(Y, P, \Lambda) \to H_4(Y, K, \Lambda)$. Because $H_4(Y, \Lambda) \to H_4(Y, P, \Lambda)$ is an isomorphism (here we use our assumption $H_4(P, \Lambda) \cong 0$ and Lemma 4.6), the image of $H_4(Y, \Lambda)$ in $H_4(Y, K, \Lambda)\equiv \Lambda \oplus \Lambda$ is generated by $(0,1)$. So we can write the diagram 
$$
\CD
@. H_4(Y, \Lambda) \\
@. @VVV \\
\Lambda  @>>> \Lambda \oplus \Lambda @>>> (\Lambda \oplus \Lambda)/\Lambda (1,1) \\
@| @|  \\
\Lambda @>>> \Lambda \oplus \Lambda @>>> \Lambda
\endCD
$$
The map $\bar j_{*}$ corresponds to $\Lambda \to \Lambda \oplus \Lambda$ defined by $1\to (1,1)$. Hence the map $H_4(Y, \Lambda) \to H_4(Y, X, \Lambda)$ corresponds to the isomorphism $\Lambda \to (\Lambda \oplus \Lambda)/\Lambda (1,1)$ defined by $1\to [(0,1)]$, the class of $(0,1)$ in the quotient. Therefore we have an isomorphism $H_3(X, \Lambda) \to H_3(Y, \Lambda)$.
\qed\enddemo

\proclaim{Lemma 4.8}  The map $\cap [X] : H^1(X, \Lambda) \to H_3(X, \Lambda)$ is an isomorphism.
\endproclaim

\demo{Proof} This follows from the diagram
$$
\CD
H^1(X, \Lambda) @> \cap [X] >> H_3(X, \Lambda) \\
@A j^{*} A \cong A @V \cong V j_{*} V \\
H^1(Y, \Lambda) @> \cap j_{*} [X] >> H_3(Y, \Lambda) \\
@V i^{*} V \cong V @A \cong A i_{*} A \\
H^1(P, \Lambda) @> \cap [P] > \cong > H_3(P, \Lambda) 
\endCD
$$
and $\cap j_{*} [X] = \cap i_{*} [P] : H^1(Y, \Lambda)\to H_3(Y,\Lambda)$.
\qed\enddemo

\medskip

(III) {\sl Case } $p=3$

\medskip
We look now to the case $\cap [X] : H^3(X, \Lambda) \to H_1(X, \Lambda) \cong 0$, i.e., we have to show that $H^3(X, \Lambda) \cong 0$. Note that the sequence
$$
\CD
0 @>>> H^3(K, \Lambda) @>>> H^4(P, K, \Lambda) @>>> H^4(P, \Lambda) @>>> 0
\endCD
$$
is exact. Since $H^4(P, \Lambda) \cong H_0 (P, \Lambda) \cong \Bbb Z$, this sequence coincides with 
$$
\CD
0 @>>> I(\Lambda) @>>> \Lambda @> \epsilon >> \Bbb Z @>>> 0
\endCD
$$ 
where $\epsilon$  is the augmentation, and $I(\Lambda) = \operatorname{Ker} \epsilon$. Let us consider the diagram
$$
\CD
H^3(K, \Lambda) @>>> H^4(Y, P, \Lambda) @> \bar \varphi^{*} > \cong > H^4(D^4, \Bbb S^3, \Lambda)\\
@| @AAA \\
H^3(K, \Lambda) @>>> H^4(Y, K, \Lambda) @> \bar j^{*} >> H^4(X, K, \Lambda)\\
@| @VV \bar i_{*} V \\
H^3(K, \Lambda) @>>> H^4(P, K, \Lambda) @> \bar \varphi_{1}^{*} > \cong > H^4(D^4, \Bbb S^3, \Lambda)
\endCD
$$
The two vertical maps split $H^4(Y, K, \Lambda) \cong \Lambda \oplus \Lambda$ such that
$$
\bar i^{*} : H^4(Y, K, \Lambda)\cong \Lambda \oplus \Lambda \to \Lambda \cong H^4(P, K, \Lambda)
$$
projects onto the first component, and $H^4(Y, K, \Lambda) \to H^4(Y, P, \Lambda)\cong \Lambda$ projects onto the second component. Because the composition $H^3(K, \Lambda) \to H^4(Y, P, \Lambda)$ is the zero map, we can identify the image of $H^3(K, \Lambda)\to H^4(Y, K, \Lambda)$ with $(I(\Lambda), 0)\subset \Lambda \oplus \Lambda$. The map $\bar j^{*}$ is the sum $\Lambda \oplus \Lambda \to \Lambda$  since the generator of $H^4(X, K, \Lambda)\cong \Lambda$ maps under
$$
\bar \varphi_1^{*} + \bar \varphi^{*} : H^4(X, K, \Lambda)\to H^4(D^4, \Bbb S^3, \Lambda)\cong \Lambda
$$
to a generator. Hence the image of
$$
\CD
H^3(K, \Lambda) @>>> H^4(Y, K, \Lambda) @> \bar j^{*} >> H^4(X, K, \Lambda)
\endCD
$$ 
is $I(\Lambda)\subset \Lambda$, i.e., $H^3(K, \Lambda) \to H^4(X, K, \Lambda)$ is injective. The long exact sequence of the pair $(X, K)$ implies $H^3(X, \Lambda) \cong 0$.
\medskip

(IV) {\sl Case} $p=4$

\medskip
{\sl Remark}. The last argument implies also $H^4(X, \Lambda) \cong \Lambda /I(\Lambda) \cong \Bbb Z$. We have proved the first part of the following

\proclaim{Lemma 4.9}  $H^3(X, \Lambda)\cong 0$, $H^4(X, \Lambda)\cong \Bbb Z$, and $\cap [X] : H^4(X, \Lambda)\to H_{0}(X, \Lambda)$ is an isomorphism. 
\endproclaim

\demo{Proof} The second part follows from the well--known property of cap products indicated in the diagram
$$
\CD
\Bbb Z \cong H^4(X, \Lambda) @> \cap [X] >> H_{0}(X, \Lambda)\cong \Bbb Z  \\
@A \epsilon AA  @AA \epsilon A \\
\Lambda \cong H^4(X, K, \Lambda) = \operatorname{Hom}_{\Lambda} (C_{4}(\widetilde X, \widetilde K), \Lambda) @> A > \cong > C_{0}(\widetilde X) \cong \Lambda
\endCD
$$
Here $A(\alpha) = \alpha (1)$, being $1\in C_4(\widetilde X, \widetilde K)$ the generator. Observe that $H_{0}(X, \Lambda) = C_{0}(\widetilde X)/ \partial_1 C_1 (\widetilde X)$, so $\epsilon$ corresponds to the canonical map $C_{0}(\widetilde X) \to C_{0}(\widetilde X)/ \partial_1 C_1(\widetilde X)$ (we may assume that $X$ has one $0$--cell).
\qed\enddemo

\medskip

(V) {\sl Case} $p=2$.

\medskip

Recall the $2$--stage Postnikov system for $P$
$$
\CD
P @> f >> D \\
@| @VV \chi V \\
P @> p >> B\pi
\endCD
$$
Let $f_{0} = f|_{K}$. Given any $\psi : \Bbb S^3 \to K$, a canonical map $g : K\cup_{\psi} D^4 \to D$ can be constructed as follows. Let $H : \Bbb S^3 \times I \to D$ be the zero homotopy of the composition $f_{0} \circ \psi : \Bbb S^3 \to D$. It factors over 
$$
\CD
D^4 = (\Bbb S^3 \times I)/\Bbb S^3 \times \{1\} @> \hat H >> D.
\endCD
$$
Then $g = f_{0} \cup \hat H : K \cup_{\psi} D^4 \to D$. Since $\pi_q (D)\cong 0$ for $q\ge 3$, $g$ is unique up to homotopy. In our case we have $\psi = \varphi_2 = \varphi_1 + \varphi$ with $[\varphi] \in \Gamma(\pi_2)$, where $\varphi : \Bbb S^3 \to K^{(2)}$, i.e., we need the zero homotopy of the composition
$$
\CD
\Bbb S^3 @>>> \Bbb S^3 \lor \Bbb S^3 @> \varphi_1 \lor \varphi >> K \lor K^{(2)} @> f_{0} \lor f_{0} >> D\lor D @>>> D.
\endCD
$$
We take the wedge of the zero homotopies $H : \Bbb S^3 \times I\to D$ for $f_{0}\circ \varphi_1$ and
$H_{0} : \Bbb S^3 \times I\to D$ for $f_{0}\circ \varphi$. This gives us the following
\medskip

\proclaim{Lemma 4.10}  Let $g_{0} = f_{0} \cup \hat H_{0} : K^{(2)} \cup_{\varphi} D^4 \to D$ denote the canonical extension  and $\theta \in H_4(K^{(2)} \cup_{\varphi} D^4, \Bbb Z)$ the canonical generator. Then we have
$$
g_{*} [X] = f_{*} [P] + (g_{0})_{*} (\theta).
$$
\endproclaim

\proclaim{Corollary 4.11}  $(g_{0})_{*}(\theta) \in \operatorname{Ker} G \subset H_4(D, \Bbb Z)$. In particular,
$$
\cap f_{*} [P] = \cap g_{*} [X] : H^2(D, \Lambda) \to H_2(D, \Lambda)
$$
that is, the map $\cap [X] : H^2 (X, \Lambda) \to H_2(X, \Lambda)$ is an isomorphism.
\endproclaim

\demo{Proof} The above spectral sequence applied to $K^{(2)} \cup_{\varphi} D^4$ gives
$$
\CD
0 @>>> \Bbb Z \otimes_{\Lambda} H_4(K^{(2)} \cup_{\varphi} D^4, \Lambda) @>>> H_4(K^{(2)} \cup_{\varphi} D^4, \Bbb Z) \\
@.  @>>> H_2(B\pi, H_2(K^{(2)} \cup_{\varphi} D^4, \Lambda)) @>>> 0.
\endCD
$$
The first map is an isomorphism, so $H_2(B\pi, H_2(K^{(2)} \cup_{\varphi} D^4, \Lambda))\cong 0$. Compared with the exact sequence for $D$
$$
\minCDarrowwidth{12pt}
\CD
0 @>>> \Bbb Z \otimes_{\Lambda} H_4(K^{(2)} \cup_{\varphi} D^4, \Lambda) @>>> H_4(K^{(2)} \cup_{\varphi} D^4, \Bbb Z)
@>>> 0 \\
@. @. @VV (g_{0})_{*} V @VVV \\
@. @. G: H_4(D, \Bbb Z) @>>> \operatorname{Hom}_{\Lambda} (H^2(D, \Lambda), H_2(D, \Lambda))
\endCD
$$
gives the result.
\qed\enddemo

\medskip

\proclaim{Theorem 4.12}  Suppose $B\pi$ is homotopy equivalent to a $2$--dimensional complex. Let $\pi_2 = H^2(B\pi, \Lambda)$. Then, if we fix one model $P$, we obtain all models by the above construction.  
\endproclaim

\demo{Proof} Fixing $P$, we constructed for any $[\varphi ]\in \pi_2$ a strongly minimal model.
Conversely, let $X = K \cup_{\psi} D^4$ be a minimal model, where $\psi : \Bbb S^3 \to K$ is the attaching map. The map $f : X\to D$ into the 2-stage Postnikov space $D$ is given by the zero homotopy of 
$$
\CD
\Bbb S^3 @> \psi >> K @> f_{0} >> D
\endCD
$$
that is
$$
\CD
\Bbb S^3 \times I @> H >> D\\
@VVV @AA \hat H A \\
D^4 = (\Bbb S^3 \times I)/\Bbb S^3 \times \{ 1\} @= D^4
\endCD
$$
with $f = f_{0} \cup \hat H$. Let us consider $\hat H : (D^4, \Bbb S^3) \to (D, K)$ and let 
$$
\bar \psi : (D^4, \Bbb S^3) \to (X, K)
$$ 
be the top cell. The diagram
$$
\CD
H_4(X, \Bbb Z) @>> \cong > H_4(X, K, \Bbb Z) @< \bar \psi < \cong < H_4(D^4, \Bbb S^3, \Bbb Z) \\
@V f_{*} VV @. @| \\
H_4(D, \Bbb Z) @>>> H_4(D, K, \Bbb Z) @< \hat H_{*} << H_4(D^4, \Bbb S^3, \Bbb Z) 
\endCD
$$
shows that $f_{*} [X]$ depends only on $\psi \otimes_{\Lambda} 1 \in \pi_3(K) \otimes_{\Lambda} \Bbb Z$. Note that $H_4(D, \Bbb Z) \to H_4(D,K, \Bbb Z)$ is injective. This also demonstrates that the above construction only depends on $\xi$, not on the choice of $[\varphi] \in \Gamma(\pi_2)$ with $[\varphi] \otimes_{\Lambda} 1 =\xi$.

It remains to show that any minimal model $X'$ is homotopy equivalent to some model $X$ obtained by the above construction. Write 
$$
\CD
X' = K' \cup_{\psi} D^4 @> f' >> D' 
\endCD
$$
where $D'$ is the $2$--stage Postnikov space, $K'$ is a $3$--dimensional complex and $\psi : \Bbb S^3 \to K'$ is the attaching map. Recall our standard model 
$$
\CD
P = K \cup_{\varphi_1} D^4 @> f >> D
\endCD
$$
In Section 3 we constructed a homotopy equivalence $h : D'\to D$ sending 
$K'\to K$. Lemma 3.3 implies 
$$
h_{*} f'_{*} [X'] - f_{*} [P] \in \operatorname{Ker} G = \Gamma(\pi_2)\otimes_{\Lambda} \Bbb Z
$$
By Lemma 4.1 of Section 4 choose $[\varphi] \in \Gamma(\pi_2)$ so that $[\varphi] \otimes_{\Lambda} 1 = h_{*} f'_{*} [X'] - f_{*} [P]$, and $\varphi : \Bbb S^3 \to K^{(2)} \subset K$. As in Part V of Section 4 we build $X = K\cup_{\varphi_2} D^4$, with $\varphi_2 = \varphi_1 + \varphi$, and $g : X\to D$. Let $g_{0} : K^{(2)} \cup_{\varphi} D^4 \to D$ be the canonically defined map from the zero homotopy of $\Bbb S^3 \to K^{(2)} \to D$. Then we have (use Lemma 4.10) $g_{*} [X] = f_{*} [P] + (g_{0})_{*} (\theta)$, where $\theta \in H_4(K^{(2)} \cup_{\varphi} D^4, \Bbb Z)$ is a generator. But $(g_{0})_{*} (\theta) = h_{*} f'_{*} [X'] - f_{*} [P]$ as it can be seen from the diagram
$$
\minCDarrowwidth{12pt}
\CD
0 @>>> H_4(K^{(2)} \cup_{\varphi} D^4, \Lambda) \otimes_{\Lambda} \Bbb Z @>> \cong > H_4(K^{(2)} \cup_{\varphi} D^4, \Bbb Z) \\
@. @V (g_{0})_{*} \otimes_{\Lambda} 1 VV @VV (g_{0})_{*} V \\
0 @>>> \Gamma(\pi_2) \otimes_{\Lambda} \Bbb Z = H_4(D, \Lambda) \otimes_{\Lambda} \Bbb Z @>>> H_4(D, \Bbb Z) @>>> H_2(B\pi, H_2(D, \Lambda)) @>>> 0.
\endCD
$$ 
Therefore $g_{*} [X] = h_{*} f'_{*} [X']$. By Proposition 3.2 and the proof of Theorem 1.3 (where we have to use that $\pi_2$ is not a torsion group) we obtain a homotopy equivalence $X'\to X$.
\qed\enddemo

\medskip

\subhead 5. Non--uniqueness of strongly minimal models: examples
\endsubhead
\bigskip

In Section 4 we constructed minimal models for all elements of $\Gamma(\pi_2)$. In this section we address the question of uniqueness up to homotopy equivalence. Recall that for two models $X$ and $X'$ we have a homotopy equivalence between the $2$--stage Postnikov systems (assuming that the first $k$--invariants are compatible). It is deduced from Diagram (2) in Section 3, i.e., we have a diagram 
$$
\CD
X @> f >> D @> h >> D'  @< f' << X' \\
@. @VVV @VVV \\
@. B\pi_1 @>>> B\pi_1 .
\endCD
$$
If $X = K \cup_{\varphi} D^4$ and $X' = K' \cup_{\psi} D^4$, then $D$ and $D'$ are constructed from the $3$--complexes $K$ and $K'$, respectively, by adjoining cells of dimension greater or equal to $4$. Proposition 3.2 defines an obstruction to extend the restriction $h^{(3)} : K\to K'$ to a homotopy equivalence $X\to X'$. Also, if this obstruction does not vanish, it could be that $X$ is homotopy equivalent to $X'$. We use $h$ to identify $D\to B\pi_1$ with $D'\to B\pi_1$. All this makes sense if $B\pi_1$ is an aspherical $2$--complex. From now on we shall consider only Baumslag--Solitar groups $B(k)$, $k\neq 0$, and aspherical surface fundamental groups. For any such model $X$ we obtain $H_3(X, \Lambda) \cong H^1(X, \Lambda) \cong H^1(B\pi, \Lambda) \cong 0 $ by Lemma 6.2 of [5] (here $\pi =\pi_1$, as usual). Since $H_4(X,\Lambda)\cong 0$, we get an isomorphism from $H_4(X,K, \Lambda)$ onto $H_3(K, \Lambda)$, i.e., $H_3(K, \Lambda)\cong \Lambda$. Furthermore, the canonical generator of $H_4(X, K, \Lambda)$, given by the attaching map $\varphi$ defines a generator of $H_3(K, \Lambda)$ and a splitting $s_X : H_3(K, \Lambda) \to \pi_3(K)$ of the Whitehead sequence given by the diagram:
$$
\CD
0 @>>> \Gamma (\pi_2) @> i_{*} >> \pi_3(K) @> H >> H_3(K, \Lambda) @>>> 0 \\
@. @. @AAA @AA \cong A \\
@. @.  \pi_4(X, K) @> \cong >> H_4(X, K, \Lambda). 
\endCD
$$  
Then $s_X$ defines a splitting $t_X : \pi_3(K) \to \Gamma (\pi_2)$. From the Whitehead sequence of $X$, we have an isomorphism from $\Gamma (\pi_2)$ onto $\pi_3(X)$, and $t_X$ can also be defined by the diagram
$$
\CD
\Gamma (\pi_2) @< t_X << \pi_3(K) \\
@V \cong VV @VV j_{*} V \\
\pi_3(X) @= \pi_3(X).
\endCD
$$
Conversely, $t_X$ defines $s_X$ by the well--known procedure, using the projection operator $i_{*}\circ t_X$. 
If $X = K \cup_{\varphi} D^4$ and $X' = K \cup_{\psi} D^4$ are homotopy equivalent models, there is a homotopy equivalence of pairs (see [13], Theorem 2.4)
$$
g : (K, \varphi (\Bbb S^3)) \to (K, \psi (\Bbb S^3))
$$
inducing the diagrams
$$
\CD
0 @>>> \Gamma (\pi_2) @> i_{*} >> \pi_3(K) @> j_{*} >> \pi_3(X)\\
@. @V g_{*} VV @V g_{*} VV @V g_{*} VV \\
0 @>>> \Gamma (\pi_2) @> i'_{*} >> \pi_3(K) @> j'_{*} >> \pi_3(X')
\endCD
$$
and
$$
\CD
0 @>>> \Gamma (\pi_2) @> i_{*} >> \pi_3(K) @> H >> H_3(K, \Lambda) @>>> 0 \\
@. @V g_{*} VV @V g_{*} VV @V g_{*} VV \\
0 @>>> \Gamma (\pi_2) @> i'_{*} >> \pi_3(K) @> H >> H_3(K, \Lambda) @>>> 0.
\endCD
$$
Hence all splittings $t_X$, $t_{X'}$, $s_X$ and $s_{X'}$ commute with the induced homomorphisms $g_{*}$. In the following we fix one model $X = K \cup_{\varphi} D^4$. We are going to construct models  $X' = K \cup_{\psi} D^4$ which are not homotopy equivalent to $X$. Let us denote by $1 \in H_3(K, \Lambda)$ the generator defined by $X$, i.e., $s_X(1) = [\varphi]$. Let $\theta : \Gamma (\pi_2) \to \Gamma (\pi_2)$ be an isomorphism. Then $\theta \circ t_X = t: \pi_3(K)\to \Gamma(\pi_2)$ is a 
splitting. It defines a splitting $s : H_3(K, \Lambda) \to \pi_3(K)$. Then $ s(1) = s_X (1) + i_{*} (a)$ for some $a \in \Gamma (\pi_2)$. As in Section 4, we construct the model $X' = K \cup_{\psi} D^4$ with $[\psi ] = s(1)$.

\medskip

\proclaim{Proposition 5.1}  If $\theta$ is not induced from an isomorphism $\pi_2 \to \pi_2$, then $X'$ is not homotopy equivalent to $X$
\endproclaim

\demo{Proof} Any homotopy equivalence $g: X\to X'$ induces
$$
\CD
0 @>>> \Gamma (\pi_2) @>>> \pi_3(K) @>>> \pi_3(X)\\
@. @V g_{*} VV @V g_{*} VV @V g_{*} VV \\
0 @>>> \Gamma (\pi_2) @>>> \pi_3(K) @>>> \pi_3(X').
\endCD
$$
However, $g_{*} : \Gamma (\pi_2) \to \Gamma (\pi_2)$ is never $\theta$.
\qed\enddemo

\noindent
{\sl Examples}. Let $X = F\times \Bbb S^2$, where $F$ is a closed oriented aspherical surface. Then $\pi_2(X) \cong \Bbb Z$, $\Gamma (\pi_2) \cong \Bbb Z$ and $- \operatorname{Id} : \Gamma (\pi_2) \to \Gamma (\pi_2)$ is not induced from an isomorphism $\pi_2 \to \pi_2$. This follows easily from the $\Gamma$--functor property. There are inclusions $\pi_2 \to \Gamma (\pi_2)$ and $\Gamma (\pi_2) \to \pi_2 \otimes \pi_2$ (because $\pi_2$ is free abelian) such that the composition $\pi_2 \to \Gamma (\pi_2) \to \pi_2 \otimes \pi_2$ sends $x$ to $x\otimes x$. In the case, $\pi = B(k)$, $\pi_2$ is free abelian (see [5], Lemma 6.2 V); one obtains such $\theta$ in this case, too. On the other hand, if $\theta$ is induced from an isomorphism $\beta : \pi_2 \to \pi_2$, one needs more to construct a homotopy equivalence. By [15], Theorem 3, one gets a map $g : K\to K$, but the induced maps $g_{*}$ not necessarily commute with the splittings $s_X$ and $s_{X'}$.   

\medskip

\noindent
{\bf Supplement to the aspherical surface case}. In the example $F\times \Bbb S^2$ there are two models, namely, $F\times \Bbb S^2$ and the non--trivial $\Bbb S^2$--bundle $E\to F$ with second Stiefel--Whitney class $\neq 0$ (see, for example, [3], Appendix). Here it is also convenient to consider the map
$$
F_{\Bbb Z} : H_4(D, \Bbb Z) \to \operatorname{Hom}_{\Bbb Z} (H^2(D, \Bbb Z) \otimes H^2(D, \Bbb Z), \Bbb Z)
$$
given by
$$
F_{\Bbb Z} (x) (u\otimes v) : = x \cap (u \cup v)
$$
where $D = F\times \Bbb CP^{\infty}$. Then $F_{\Bbb Z}$ is injective. If $f_0 : F\times \Bbb S^2 \to D$ and $f_1 : E\to D$ are the Postnikov maps, then $F_{\Bbb Z}(f_{0 *}[F\times \Bbb S^2])$ and $F_{\Bbb Z}(f_{1 *}[E])$ are the integral intersection forms of $F\times \Bbb S^2$ and $E$, respectively. Moreover, these forms are given by the matrices
$$
\pmatrix 0 & 1 \\
1 & 0 \endpmatrix \qquad \qquad 
\pmatrix 0 & 1 \\
1 & 1 \endpmatrix
$$
respectively (see [3]). It is shown in [9], Section 5, that $F\times \Bbb S^2$ and $E$ are the only models up to homotopy equivalence.
\bigskip

\subhead 6. Final remarks
\endsubhead
\bigskip

The following map was defined in [2] 
$$
F : H_4(D, \Bbb Z) \to \operatorname{Hom}_{\Lambda- \Lambda} (H^2(D, \Lambda) \otimes_{\Bbb Z} \overline{H}^2(D, \Lambda), \Lambda)
$$
to classify Poincar\'e $4$--complexes $X$, where $D\to B\pi$ is a $2$--stage Postnikov system for $X$. Here $H^2(D, \Lambda) \otimes_{\Bbb Z} \overline{H}^2(D, \Lambda)$ carries the obvious $\Lambda$--bimodule structure. It was  proved therein that $F$ is injective for free non-abelian groups $\pi$. The maps $F$ and $G$ are related by the following diagram
$$
\CD
H_4(D, \Bbb Z) @> G >> \operatorname{Hom}_{\Lambda- \Lambda} (H^2(D, \Lambda), \overline{H}_2(D, \Lambda)) \\
@| @VV H V \\
H_4(D, \Bbb Z) @> F >> \operatorname{Hom}_{\Lambda- \Lambda} (H^2(D, \Lambda) \otimes_{\Bbb Z} \overline{H}^2(D, \Lambda), \Lambda)
\endCD
$$
where $H(\varphi ) (u\otimes v) = \overline{\hat u (\varphi (v))}$, and $\hat u$ is the image of $u$ under 
$$
H^2(D, \Lambda)\to \operatorname{Hom}_{\Lambda}(H_2(D,\Lambda), \Lambda).
$$
Obviously, $G$ is injective if $F$ is injective. If $f: X\to D$ is a map such that $f_{*} : \pi_q(X)\to \pi_q(D)$ is an isomorphism for $q= 1, 2$, then $F(f_{*} [X]) \circ (f^{*} \otimes f^{*})$ is the equivariant intersection form on $X$, and $f_{*} G (f_{*} [X]) f^{*} : H^2(X, \Lambda) \to \overline{H}_2 (X, \Lambda)$ is the Poincar\'e duality isomorphism. It is convenient to denote $F(f_{*} [X])$ the "intersection type" and $G(f_{*} [X])$ the "Poincar\'e duality type" of $X$. The Poincar\'e duality type determines the intersection type. In this sense it is a stronger "invariant". For $\Bbb S^2$--bundles over aspherical $2$--surfaces all intersection types vanish, whereas the Poincar\'e types are non-trivial.

\bigskip
{\bf Acknowledgements}.
Work performed under the auspices of the  G.N.S.A.G.A. of the
C.N.R. (National Research Council)  of Italy  and partially
supported by  the MIUR (Ministero per la Ricerca Scientifica e Tecnologica) of Italy,  and by the  Slovenia Research Agency grants  No. P1--0292--0101 and J1--2057--0101.

\medskip

\noindent
Authors' addresses:

\noindent
Alberto Cavicchioli, Dipartimento di Matematica, Universit\`a di Modena
e  Reggio Emilia, Via Campi 213/B, 41100 Modena, Italy. 
E-mail: alberto.cavicchioli\@unimore.it

\medskip

\noindent
Friedrich Hegenbarth, Dipartimento di Matematica, Universit\`a di
Milano, Via Saldini 50, 20133 Milano, Italy. 
E-mail: friedrich.hegenbarth\@unimi.it

\medskip

\noindent
Du\u san Repov\u s, Faculty of Mathematics and Physics, 
University of Ljubljana, P.O.B. 2964, Ljubljana, Slovenia 1001.
E-mail: dusan.repovs\@guest.arnes.si 

\bigskip

\enddocument        

\Refs

\ref 
\no 1
\by   Baues, H. J.:
Combinatorial Homotopy and $4$--Dimensional Complexes.
 Walter de Gruyter.  Berlin --New York, 1991 
\endref



\ref
\no 2
\by  Cavicchioli, A.,  Hegenbarth, F.:
On $4$--manifolds with free fundamental groups.
Forum  Math. 6, 415--429 (1994).
\endref


\ref
\no 3  
\by Cavicchioli, A., Hegenbarth, F.,  Repov\u s, D.:
Four--manifolds with surface fundamental groups.
Trans. Amer. Math. Soc. 349, 4007--4019 (1997).
\endref


\ref
\no 4
\by Hambleton, I., Kreck, M.:
On the classification of topological 
$4$--manifolds with finite fundamental group.
Math. Ann.  280, 85--104 (1988).
\endref

\ref
\no 5
\by I. Hambleton, I., Kreck, M., Teichner, P.:
Topological $4$--manifolds with geometrically $2$--dimensional fundamental group.
J. Topology Anal.   1, 123--151 (2009).
\endref

\ref
\no 6
\by  Hegenbarth, F., Repov\u s, D., Spaggiari, F.:
Connected sums of $4$--manifolds.
Topology Appl.   146-147,  209--225 (2005).
\endref

\ref
\no 7
\by Hegenbarth, F., Piccarreta, S.:
 On Poincar\'e $4$--complexes with free fundamental groups.
Hiroshima Math. J.   32, 145--154 (2002).
\endref

\ref
\no 8
\by   Hillman, J.A.:
$PD_4$--complexes with strongly minimal models.
 Topology  Appl.  153, 2413--2424 (2006).
\endref


\ref
\no 9
\by Hillman, J.A.:
Strongly minimal $PD_4$--complexes.
Topology Appl.  156, 1565--1577 (2009).
\endref

\ref
\no 10
\by Kaplansky, I.R.:
Fields and Rings.
Chicago Univ. Press. Chicago--London, 1969.
\endref


\ref
\no 11
\by  Mac Lane, S.,  Whitehead, J.H.C.:
 On the 3--type of a complex.
Proc. Nat. Acad. Sci.  36, 41--48 (1950).
\endref

\ref
\no 12
\by Pamuk, M.:
Homotopy self--equivalences of $4$--manifolds with $\operatorname{PD}_2$--fundamental group. 
Topology Appl. 156, 1445--1458 (2009).
\endref


\ref
\no 13
\by  Wall, C.T.C.:
Poincar\'e complexes.
Ann. of Math.  86, 213--245 (1967).
\endref


\ref
\no  14
\by  Wall, C.T.C.:
Surgery on Compact Manifolds.
Academic Press.
London -- New York
1970
\endref

\ref
\no 15
\by  Whitehead, J.H.C.:
On a certain exact sequence.
Ann. of Math. 52 (2), 51--110 (1950).
\endref
\endRefs
\bye